\documentclass[12pt,a4paper]{amsart}
\usepackage{amsmath}
\usepackage{hyperref}
\usepackage{graphicx}
\usepackage{ifpdf}
\ifpdf
\fi

\newcommand{\incg}[2][.5in]{\setbox5=\hbox{\;\includegraphics[height=#1]{#2}\;}%
\dimen1=-#1\divide\dimen1 by 2\raise\dimen1\box5}
\newcommand{\vc}[1]{\mathbf{#1}}

\newcommand{\bR}{\mathbb{R}}
\newcommand{\bC}{\mathbb{C}}
\newcommand{\bH}{\mathbb{H}}
\newcommand{\bO}{\mathbb{O}}
\newcommand{\bK}{\mathbb{K}}

\newcommand{\fe}{\mathfrak{e}}
\newcommand{\ff}{\mathfrak{f}}
\newcommand{\fg}{\mathfrak{g}}

\newcommand{\fsl}{\mathfrak{sl}}

\newcommand{\fder}{\mathfrak{der}}

\newcommand{\SL}{\mathrm{SL}}

\newcommand{\Spin}{\mathrm{Spin}}

\newcommand{\End}{\mathrm{End}}
\newcommand{\tr}{\mathrm{tr}}
\newcommand{\Cliff}{\mathrm{Cliff}}

\newcommand{\upabit}[1]{\setbox5=\hbox{\;\ensuremath{#1}\;}\raise7pt\box5}

\begin{document}
\title{Hurwitz' theorem on composition algebras}
\author{Bruce W. Westbury}
\address{%
Mathematics Institute\\
University of warwick\\
Coventry CV4 7AL}
\email{Bruce.Westbury@warwick.ac.uk}
\date{2009}
\maketitle

\section{Introduction}
In this article we describe several results based on the paper \cite{hurwitz}
and which we will refer to as Hurwitz' theorem. There are several related results:
the classification of real normed division algebras, the classification of complex composition
algebras and the classification of real composition algebras.
Our interest is in real composition algebras. There are seven real composition algebras.

This article was originally conceived of as an exposition of the proof in \cite{Diploma} (available as
\href{http://www.mathematik.uni-bielefeld.de/~rost/data/boos.pdf}{Diploma}
) and \cite{MR1397790} 
and the intention was to promote the use of diagrams. This result can also be found in \cite[Chapter 16]{MR2418111} although it is less explicit. This is still the core of this article.
However I then added a prequel giving some background on composition algebras and vector product
algebras. In particular I give brief accounts of two other proofs of Hurwitz' theorem. Both proofs
can be found with more detail in \cite[\S 2]{MR1886087}. I then added a sequel based on 
\cite[Volume 1]{MR1701618}. The relevant Chapters are Chapter 6 in Notes on Spinors and Chapters 1 \& 2 in
Supersolutions. This is also explained in more detail in \cite{YangMills}.
The aim here is to explain an algebraic property of the Poincare supergroup and
Minkowski superspace and to show that this is satisfied only for Minkowski space of dimension
$3$, $4$, $6$ and $10$. As a stepping stone to this I have also included a section with an
unconventional view on triality. 

In these notes we work over a field $F$ whose characteristic is not two and for \eqref{pol}
we require that the characterisric is not two or three.
\section{Algebras}
In this section we discuss composition algebras and vector product algebras.
The first result is that these have the same classification. Then we also give
a classification of composition algebras based on Cayley-Dickson doubling and a 
classification of vector product algebras based on the structure theory of Clifford algebras.
For these classifications we work over the field of complex numbers.

An {\em algebra} is a vector space $A$ with a bilinear multiplication or equivalently
a linear map $\mu\colon A\otimes A\rightarrow A$. An inner product on the vector space $A$
is {\em associative} if
\[ \left\langle x,y\times z\right\rangle = \left\langle x\times y, z\right\rangle \]
for all $x,y,z\in A$.

If $A$ has a unit $1$ then the {\em imaginary part} is
\[ \{ a\in A | \left\langle a,1\right\rangle =0 \} \]
The imaginary part has a multiplication given by
\[ u\times v = uv - \left\langle u,v\right\rangle.1 \]

Let $A$ be an algebra. A {\em derivation} of $A$ is a linear map $\partial\colon A\rightarrow A$
that satisfies the Leibniz rule
\[ \partial(ab)=\partial(a)b + a\partial(b) \]
for all $a,b\in A$. Then $\fder(A)$ is the Lie algebra of derivations of $A$.

If $A$ has an associative inner product then we define $D\colon A\otimes A\rightarrow \fder(A)$ by
\[ \left\langle  \partial , D(a,b)\right\rangle = \left\langle \partial(a),b\right\rangle \]
for all $a,b\in A$, $\partial\in\fder(A)$.

\subsection{Composition algebras}
A {\em composition algebra} has a bilinear multiplication, a unit, a symmetric inner product and an anti-involution called conjugation. The first condition is
\[ \qquad \overline{a}=-a+2\left\langle a,1\right\rangle.1 \]
In particular the conjugation is determined by the functional $a\mapsto \left\langle a,1\right\rangle$.
The second condition is 
\[ a\overline{b}+b\overline{a}=\left\langle a,b\right\rangle.2 \]
In particular the inner product is determined by the conjugation. This condition is usually
written as $a\overline{a} = |a|.1$ whereas we have preferred the polarised version.
The final condition is the key condition. This is usually written as
\[ \qquad |ab|=|a|.|b|  \]
The polarisation of this condition is
\[ 2\left\langle  a,b\right\rangle \left\langle  c,d\right\rangle
= \left\langle  ac,bd\right\rangle + \left\langle  ad,bc\right\rangle\]

Let $A$ be a composition algebra. Then we define the Cayley-Dickson double of $A$, $D(A)$.
The underlying vector space is $A\oplus A$. The multiplication is given by
\[ (a,b)(c,d)=(ac-d\overline{b},\overline{a}d+cb) \]
The unit is $(1,0)$.
The conjugation is given by
\[ \overline{(a,b)}=(\overline{a},-b) \]
The inner product is given by
\[ \left\langle  (a,b),(c,d)\right\rangle =
\left\langle  a,c\right\rangle + \left\langle  b,d\right\rangle \]

This is the structure needed for a composition algebra but the conditions may not be met.
The precise result is that $D(A)$ is a composition algebra if and only if multiplication
in $A$ is associative.

The basic example of a composition algebra over a field $K$ is the field $K$ itself.
The multiplication is the field multiplication and this has a unit. The conjugation map
is the identity map. The inner product is then determined and is given by
$\left\langle  a,b\right\rangle = ab$.

Take the field to be the field of real numbers, $\bR$. Then starting with $\bR$ and applying
the Cayley-Dickson doubling, the algebras that are constructed are the real numbers $\bR$,
the complex numbers $\bC$, the quaternions $\bH$ and the octonions $\bO$. The multiplication
in $\bO$ is not associative and the Cayley-Dickson double of $\bO$ is not a composition algebra.

Taking the field to be the field of complex numbers, $\bC$, gives the complex versions of these
algebras. The complex version of $\bR$ is $\bC$, the complex version of $\bC$ is $\bC\oplus \bC$
and the complex version of $\bH$ is the algebra of $2\times 2$ matrices with entries in $\bC$.

Next we give an outline of the classification of composition algebras over the complex numbers.
The result is that the four examples we have constructed using Cayley-Dickson doubling are the
only possibilities.

The proof of this result is based on the observation that if $A$ is a composition algebra
and $A_0$ is a proper subalgebra then $D(A_0)$ is a subalgebra of $A$. To see this choose
$\imath\in A$ such that $\left\langle  a,\imath\right\rangle = 0$ for all $a\in A_0$
and such that $\left\langle  \imath,\imath\right\rangle = 1$. Then we have an inclusion
 $D(A_0)\rightarrow A$ given by $(a,b)\mapsto a+\imath b$ for all $a,b\in A_0$.

Then the result follows from this. Let $A$ be a composition algebra. Then by the previous 
observation $\bO$ is not a proper subalgebra of $A$. Also $\bC$ is a subalgebra of $A$
by taking scalar mutliples of the unit. Hence it follows from the previous observation that
$A$ is isomorphic to one of the four composition algebras constructed by Cayley-Dickson 
doubling. 

There are seven real composition algebras. Four of these are obtained by Cayley-Dickson
doubling from $\bR$. This gives the real numbers $\bR$, the complex numbers $\bC$, the
quaternions $\bH$ and the octonions $\bO$. The inner product is definite and so these are
division algebras. The other three real composition algebras have indefinite inner product
and have zero divisors. These are known as the split compostion algebras. The three split
real composition algebras are $\bR\oplus \bR$, $M_2(\bR)$, and the split octonions.
The split octonions were constructed by Zorn.

Consider $\bR^3$ with the dot product and vector product. Then the split octonions as a
vector space is the set of $2\times 2$ matrices of the form
\[ \begin{pmatrix} a & u \\ v & b \end{pmatrix} \]
The multiplication is defined by
\begin{equation*}
 \begin{pmatrix} a & u \\ v & b \end{pmatrix}
\begin{pmatrix} a^\prime & u^\prime \\ v^\prime & b^\prime \end{pmatrix} =
\begin{pmatrix}
 aa^\prime +u.v^\prime & au^\prime+b^\prime u+v\times v^\prime \\
a^\prime v+bv^\prime - u\times u^\prime & bb^\prime +v.u^\prime
\end{pmatrix}
\end{equation*}
The norm is defined by
\[ \det \begin{pmatrix} a & u \\ v & b \end{pmatrix} = ab - u.v \]
and the conjugation is given by
\[ \overline{\begin{pmatrix} a & u \\ v & b \end{pmatrix}}=
\begin{pmatrix} b & -u \\ -v & a \end{pmatrix} \]

Choose a vector $u$ with $u.u=1$. Then we have an inclusion of $M_2(\bR)$ in the split
octonions by taking
\[ \begin{pmatrix} a & b \\ c & d \end{pmatrix} \mapsto
\begin{pmatrix} a & bu \\ cu & d \end{pmatrix} \]
This subalgebra containd $\bR\oplus \bR$ by taking diagonal matrices.

\subsection{Vector product algebras}
A {\em vector product algebra} has a symmetric inner product and an alternating bilinear multiplication.
The inner product is required to be associative. These are the conditions
\begin{align}
 \left\langle x,y \right\rangle &= \left\langle y,x \right\rangle \\
 x\times y &= -y\times x\\
 \left\langle x,y\times z\right\rangle &= \left\langle x\times y, z\right\rangle
\end{align}
The final condition for a vector product algebra is
\[ (x\times y)\times x = \left\langle x,x \right\rangle y - \left\langle x,y \right\rangle x \]
The polarisation of this condition is
\begin{equation}\label{vpa}
 (x\times y)\times z + (z\times y)\times x =
 2\left\langle x,z \right\rangle y
- \left\langle x,y \right\rangle z
- \left\langle z,y \right\rangle x
\end{equation}

Let $V$ be a vector product algebra. Then $L(\vc{u})$ for $\vc{u}\in V$ is defined by
\[ L(\vc{u})\colon (a,\vc{w})\mapsto (-\vc{u}.\vc{w},a\vc{u}+\vc{u}\times\vc{w}) \]
Then a direct calculation shows that
\[  L(\vc{u})^2\colon (a,\vc{w})\mapsto -(\vc{u}.\vc{u})(a,\vc{w}) \]
or equivalently
\[  L(\vc{u})L(\vc{v})+L(\vc{v})L(\vc{u})\colon (a,\vc{w})\mapsto
 -\left\langle \vc{u},\vc{v}\right\rangle (a,\vc{w}) \]
These are the defining relations for the Clifford algebra of $V$, $\mathrm{Cliff}(V)$.
This shows that $\mathrm{Cliff}(V)$ acts on $F\oplus V$. This is a strong constraint and
the classification of vector product algebras follows from the structure theory of Clifford algebras.

\subsection{Motivation}
First we discuss the relationship between composition algebras and vector product algebras.
A composition algebra determines a vector algebra by taking the imaginary part. A vector algebra
determines a composition algebra by formally adjoining a unit.
Let $V$ be a vector product algebra over a field $F$. Then we define the structure maps of a composition algebra on the vector space $F\oplus V$ as follows. This multiplication with some historical
references is given in \cite{MR723943}.

The multiplication map is given by
\[ \left(s,\vc{u})(t,\vc{v})=(st-\vc{u}.\vc{v},
s\vc{v}+t\vc{u}+\vc{u}\times\vc{v}\right) \]
The unit is $(1,\vc{0})$. The conjugation is given by
\[ \overline{(s,\vc{v})}=(s,-\vc{v}) \]
The inner product is given by
\[ \langle (s,\vc{u})|(t,\vc{v})\rangle = st+\vc{u}.\vc{v} \]

Then this structure satisfies the conditions for a composition algebra. Furthermore these are
inverse operations. Taking the imaginary part of this composition algebra recovers the vector
product algebra. Also applying this construction to the imaginary part of a composition algebra
recovers the composition algebra.

Also the composition algebra and the vector product algebra have the same derivation Lie algebra.

Next we discuss two situations in the study of the exceptional simple Lie algebras in which
vector product algebras arise. Although it may not be apparent these are closely related.
One situation is taken from \cite[(16.11)]{MR2418111}. Let $A$ be an algebra with an associative
inner product whose multiplication is anti-symmetric and whose inner product is symmetric.
Assume that $A$ considered as a representation of the derivation Lie algebra is irreducible.
Consider the endomorphism algebra of $A\otimes A$ considered as a representation
of the derivation Lie algebra of $A$.  Assume further that the dimension of this endomorphism algebra
is at most five. Then there are just two possibilities; either $A$ is a Lie algebra or $A$
is a vector product algebra.

The second situation is the Freudenthal-Tits construction
of the exceptional simple Lie algebras given in \cite{MR0214638}. A careful account of this
is given in \cite{MR2020553}.

Let $H^+$ be the algebra of
$3\times 3$ Hermitian matrices with entries in the octonions. This has a symmetric inner
product and a symmetric multiplication. These are given by
\[ \left\langle A,B\right\rangle= \tr(AB)\qquad A\circ B=\frac12 (AB+BA) \]
Here the product $AB$ is given by matrix multiplication and is not an element of $H^+$.
This is an associative inner product.

Let $H_{26}$ be the imaginary part of $H^+$. The derivation algebra of $H_{26}$ is the exceptional simple Lie algebra $\ff_4$, \cite{MR0034378}.

Now let $V$ be an algebra and
denote the derivation algebra by $\fder(V)$. Then we start with the vector space
\[ \fg = \ff_4 \oplus \fder(V) \oplus H_{26}\otimes V \]
Then the aim is to construct a Lie bracket so that $\fg$ is a simple Lie algebra.
Here $\ff_4 \oplus \fder(V)$ is already a Lie algebra and this Lie algebra acts
on $H_{26}\otimes V$. Hence to define the Lie bracket on $\fg$ it remains to construct
an anti-symmetric map of representations of $\ff_4 \oplus \fder(V)$
\[ (H_{26}\otimes V)\otimes (H_{26}\otimes V)\rightarrow \fg \]
and to check the Jacobi identity. There is also a homomorphism of representations of $\fder(V)$,
$D\colon V\otimes V\rightarrow \fder(V)$ and this is anti-symmetric. Then we extend the Lie
bracket in the only way possible
\[ \left[ A\otimes\vc{u},B\otimes\vc{v}\right] = 
 \left\langle \vc{u},\vc{v} \right\rangle D(A,B)
+\left\langle A,B \right\rangle D(\vc{u},\vc{v})
+(A*B)\otimes(\vc{u}\times\vc{v})
\]
where $A*B$ is the multiplication on $H$. This multiplication is anti-symmetric and the
Jacobi identity is satisfied if and only if $V$ is a vector product algebra.

The examples are
\begin{align*}
 \ff_4 &= \ff_4 \oplus 0 \oplus H_{26}\otimes V_0 \\
 \fe_6 &= \ff_4 \oplus 0 \oplus H_{26}\otimes V_1 \\
 \fe_7 &= \ff_4 \oplus \fsl_2 \oplus H_{26}\otimes V_3 \\
 \fe_8 &= \ff_4 \oplus \fg_2 \oplus H_{26}\otimes V_7 \\
\end{align*}

It follows from the classification of vector product algebras that these are the only examples.

For another application of vector cross product algebras see \cite{MR2492402}.
\section{Hurwitz' theorem}
First we write the tensor equations in diagram notation. The bilinear multiplication is represented
by a trivalent vertex. The conditions that the symmetric bilinear form is non-degenerate and that
the multiplication is associative imply that the tensor only depends on the isotopy class of the
planar diagram.

The main relation is \eqref{vpa} and this is rewritten in diagram notation in \eqref{fund}.
In order to translate between these notations the first step is to write \eqref{vpa}
in index notation. We introduce the multiplication tensor $\mu\colon V\otimes V\rightarrow V$
and the inner product $g$. Then writing the relation \eqref{vpa} in index notation we have
\begin{equation*}
 \mu^{xy}_r\mu^{rz}_w+\mu^{yz}_r\mu^{xr}_w=2g^{xz}\delta^y_w-g^{xy}\delta^z_w-g^{yz}\delta^x_w
\end{equation*}
Then we introduce an edge to representation $V$ and a trivalent vertex to represent $\mu$.
Then the relation can be written as
\begin{equation*}
\incg{hur.103}+\incg{hur.104}=2\incg{hur.105}-\incg{hur.106}-\incg{hur.107}
\end{equation*}
where we have omitted the labels on the internal edges which correspond to contracted indices.
Now omitting the labels on the boundary points which correspond to uncontracted indices we
obtain \eqref{fund}.

Then we write $\delta$ for the dimension of the vector algebra.
This is written in diagram notation as
\begin{equation}\label{dim}
 \incg{hur.15}=\delta
\end{equation}

Our preference is to work with planar diagrams and to regard the crossing as part of the data.
It is more usual to work with abstract graphs (with a cyclic ordering of the edges at a vertex)
and so to regard crossing as part of the formalism. From this point of view we have the following relations.
\begin{align}\label{cross}
 \incg{hur.27}&\mapsto\incg{hur.8}\\
 \incg{hur.58}&\mapsto\incg{hur.59}\\
 \incg{hur.60}&=\incg{hur.61}
\end{align}
Note these have respectively four, five and six boundary points. The first two are taken as reduction
rules. In proof of Hurwitz' theorem we only consider diagrams with four boundary points and so this
argument only uses the first reduction rule. When we study the case $\delta=7$ then we consider diagrams
with five boundary points and we use the second rule. We will not consider diagrams with six boundary
points and so we will not make use of the third relation. There is one other relation with four boundary
points which we will use. This is the relation
\begin{equation*}
\incg{hur.26}=\incg{hur.99}
\end{equation*}

Then the symmetries are written in diagram notation as follows. These are taken to be reduction rules.
\begin{equation}\label{symm}
 \incg{hur.1}\mapsto \incg{hur.2}\qquad
 \incg{hur.3}\mapsto -\incg{hur.4}
\end{equation}

Then the fundamental identity is written in diagram notation as
\begin{equation}\label{fund}
 \incg{hur.5}+\incg{hur.6}=2\incg{hur.7}-\incg{hur.8}-\incg{hur.9}
\end{equation} 

Then we want to deduce some consequences of these relations. Firstly we have the overlap
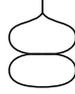
\begin{equation*}
 \incg{hur.10}
\end{equation*}
This can be simplified using either of the reduction rules \eqref{symm}.
Since the characteristic is not two this gives the reduction rule
\begin{equation}\label{tadpole}
 \incg{hur.11}=0
\end{equation}

Next we find some consequences of the fundamental identity \eqref{fund}.
The first consequence is
\begin{equation*}
 \incg{hur.45}+\incg{hur.46}=2\incg{hur.12}-\incg{hur.13}-\incg{hur.14}
\end{equation*} 
Simplifying this using \eqref{dim}, \eqref{symm}, \eqref{tadpole} gives
the reduction rule
\begin{equation}\label{bubble}
 \incg{hur.16}\mapsto(-\delta+1)\incg{hur.17}
\end{equation}

The second consequence is 
\begin{equation*}
 \incg{hur.18}+\incg{hur.19}=2\incg{hur.20}-\incg{hur.21}-\incg{hur.22}
\end{equation*} 
Simplifying this using \eqref{symm}, \eqref{tadpole}, \eqref{bubble} gives
the reduction rule
\begin{equation}\label{triangle}
 \incg{hur.23}\mapsto(\delta-4)\incg{hur.24}
\end{equation}

The third consequence is 
\begin{equation*}
 \incg{hur.25}+\incg{hur.26}=2\incg{hur.27}-\incg{hur.28}-\incg{hur.29}
\end{equation*} 
Simplifying this using \eqref{symm}, \eqref{triangle} gives
the reduction rule
\begin{equation}\label{HX}
 \incg{hur.26}\mapsto\incg{hur.5}-\incg{hur.7}-\incg{hur.9}+2\incg{hur.8}
\end{equation}

The fourth consequence is 
\begin{equation*}
 \incg{hur.33}+\incg{hur.34}=2\incg{hur.99}-\incg{hur.30}-\incg{hur.31}
\end{equation*} 
Simplifying this using \eqref{bubble}, \eqref{triangle}, \eqref{HX} gives
the reduction rule
\begin{multline}\label{sq1}
 \incg{hur.32}\mapsto (-\delta+6)\incg{hur.5}-\incg{hur.6}\\
+4\incg{hur.8}+(\delta-3)\incg{hur.9}-2\incg{hur.7}
\end{multline}

Now comes the punchline. The relation \eqref{sq1} can be rotated through a quarter of a revolution
to give the reduction rule
\begin{multline}\label{sq2}
 \incg{hur.32}\mapsto -\incg{hur.5}+(-\delta+6)\incg{hur.6}\\
+(\delta-3)\incg{hur.8}+4\incg{hur.9}-2\incg{hur.7}
\end{multline}

Then taking the difference between \eqref{sq1} and \eqref{sq2} gives the relation
\begin{equation}\label{miracle}
 (\delta -7)\left( 
\incg{hur.6}-\incg{hur.5}-\incg{hur.8}+\incg{hur.9}
\right) =0
\end{equation}

This gives two possibilities either $\delta=7$ or we have the relation
\begin{equation}\label{ass}
\incg{hur.6}-\incg{hur.5}-\incg{hur.8}+\incg{hur.9}=0
\end{equation}
The relation \eqref{ass} is equivalent to the condition that the composition algebra
associated to the vector algebra is associative. Therefore if we state this in terms of
composition algebras we have that either the dimension is 8 or the algebra is associative.

The relation \eqref{ass} has the following consequence
\begin{equation*}
\incg{hur.19}-\incg{hur.18}-\incg{hur.21}+\incg{hur.22}=0
\end{equation*}
Simplifying this gives a relation which can be compared with \cite[(16.17)]{MR2418111}
\begin{equation*}
(\delta-3)\incg{hur.4}=0
\end{equation*}
This gives two possibilities either $\delta=3$ or we have the relation
\begin{equation}\label{comm}
\incg{hur.4}=0
\end{equation}
The relation \eqref{comm} is equivalent to the condition that the composition algebra
associated to the vector algebra is commutative. 

This relation has the following consequence
\begin{equation*}
(\delta-1)\incg{hur.17}=0
\end{equation*}
This gives two possibilities either $\delta=1$ or all diagrams are $0$.

In conclusion this gives four possibilities. The completely degenerate case is
$\delta=0$ and all diagrams are $0$. The second case which is still degenerate 
is that $\delta=1$ and the following relations hold
\begin{equation*}
\incg{hur.4}=0 \qquad \incg{hur.8}=\incg{hur.9}=\incg{hur.7}
\end{equation*}
The third case is $\delta=3$ and the relation \eqref{ass} holds.
The fourth case is $\delta=7$.

\subsection{Reduction rules}
Finally we derive some further reduction rules for the case $\delta=7$.
Our aim is to deduce the finite confluent set of reduction rules in
\cite{MR1265145} and \cite{MR1403861}. An alternative approach is given in
\cite[Chapter 16]{MR2418111}. The starting point for this approach is that the 
projection onto the derivation Lie algebra is the following idempotent
in $\End(V\otimes V)$.
\begin{equation*}
\frac12\incg{hur.8}-\frac12\incg{hur.7}+\frac16\incg{hur.5}
\end{equation*}

The relation \eqref{fund} has the following consequence
\begin{equation*}
 \incg{hur.35}+\incg{hur.36}=2\incg{hur.37}-\incg{hur.38}-\incg{hur.39}
\end{equation*} 
This simplifies to give the reduction rule
\begin{equation}\label{pent1}
\incg{hur.36}\mapsto -\incg{hur.35}+2\incg{hur.37}-\incg{hur.38}-\incg{hur.39}
\end{equation}

The relation \eqref{fund} has the following further consequence
\begin{equation*}
 \incg{hur.40}+\incg{hur.41}=2\incg{hur.42}-\incg{hur.43}-\incg{hur.44}
\end{equation*} 
This simplifies to give the reduction rule
\begin{equation*}
\incg{hur.41}\mapsto -\incg{hur.40}+2\incg{hur.42}-\incg{hur.43}-\incg{hur.44}
\end{equation*} 

The fundamental relation \eqref{fund} gives the reduction rule
\begin{equation}\label{fund2}
2\incg{hur.7}\mapsto \incg{hur.5}+\incg{hur.6}
+\incg{hur.8}+\incg{hur.9}
\end{equation} 

Then we substitute $\delta=7$ in the reduction rules
\eqref{dim}, \eqref{bubble}, \eqref{triangle}
The reduction rules \eqref{sq1} and \eqref{sq2} give the reduction rule
\begin{multline}\label{sq3}
 \incg{hur.32}\mapsto -2\incg{hur.5}-2\incg{hur.6}\\
+3\incg{hur.8}+3\incg{hur.9}
\end{multline}

Now simplify the reduction rule \eqref{pent1} to obtain the reduction rule
\begin{multline*}
 \incg{hur.47}\mapsto \\
-\incg{hur.48}-\incg{hur.49}-\incg{hur.50}-\incg{hur.51}-\incg{hur.52} \\
\incg{hur.53}+\incg{hur.54}+\incg{hur.55}+\incg{hur.56}+\incg{hur.57} \\
\end{multline*}

Then we arrive at the reduction rules (with $q=1$) in \cite{MR1265145} and \cite{MR1403861}.
In particular this set of rewrite rules is confluent. Furthermore any closed diagram reduces
to a scalar by a curvature argument. Also given any $k$ the set of irreducible diagrams
with $k$ boundary points is finite; this is an application of the isoperimetric inequality.
A bijection between these finite sets and certain lattice walks is given in \cite{MR2320368}.

\section{Triality}
This is a further application of Hurwitz' theorem. 
We have seen that if $V$ is a vector product algebra then $\Cliff(V)$ acts on
$F\oplus V$. In this section and the next we promote this to higher dimensions.

Let $\bK$ be the composition algebra associated to $V$. Then we can identify $F\oplus V$
with $\bK$ and $\Cliff(V)$ with $\Cliff^{\mathrm{even}}(\bK)$. Then from the action
of $\Cliff(V)$ on $F\oplus V$ we can form an action of $\Cliff^{\mathrm{even}}(\bK)$
on $\bK$. This is usually
expressed as saying that $\bK$ is a spin representation of $\Spin(\bK)$. Next we promote this
to an action of $\Cliff(\bK)$ on $\bK\oplus \bK$ which makes this explicit.

For $a\in\bK$ define $\rho(a)\in \End(\bK\oplus \bK)$ by
\[ \rho(a)\colon (x,y)\mapsto (\overline{a}\overline{y},\overline{x}\overline{a}) \]
for all $(x,y)\in \bK\oplus \bK$. Then to check the defining relations of the Clifford
algebra we use the fact that if $x,y\in \bK$ then the subalgebra generated by
$\{1,x,y,\overline{x},\overline{y}\}$ is associative. This is equivalent to the statement
that a composition algebra is alternative. Then we calculate
\begin{align*}
 \rho(a)^2\colon (x,y)&\mapsto ( \overline{a}(ax),(ya)\overline{a}) \\
&=( (\overline{a}a)x,y(a\overline{a})) \\
&= |a| (x,y)
\end{align*}

This example has more structure. There are three vector spaces with non-degenerate inner 
products, $V_1,V_2,V_3$ and  $V_1\otimes V_2\otimes V_3\rightarrow F$.
If $(i,j,k)$ is a permutation of $(1,2,3)$
then this gives a map $V_i\otimes V_j\rightarrow V_k$. The condition is that for any
permutation $(i,j,k)$ the two maps
\[ V_i\otimes V_j\rightarrow V_k \qquad V_i\otimes V_k\rightarrow V_j \]
define an action of $\Cliff(V_i)$ on $V_j\oplus V_k$.
We will refer to this structure as a {\em triality}.

In the previous example all three vector spaces are identified with $\bK$.
The fundamental tensor is the trilinear form given by
\[ x\otimes y\otimes z\mapsto \left\langle \overline{x}\,\overline{y},\overline{z}\right\rangle 
=\left\langle xy,z \right\rangle \]
This is invariant under cyclic permutations of $(x,y,z)$.

Then we can show that this constructions gives all examples. Choose $e_1\in V_1$
and $e_2\in V_2$. Then all three spaces are identified. Take $e_3=e_1e_2$.
This gives $\bK$ with the inner product, $1\in \bK$ and the map
$x\otimes y\mapsto \overline{x}\overline{y}$. This structure is equivalent to the
structure of a composition algebra.

The diagrams for a triality are straightforward. There are three possible colours each edge.
We consider trivalent graphs coloured so that for each vertex all three edges at the vertex
have different colours. For the relations we take the one relation
\begin{equation}\label{cliff}
\incg{hur.62}+\incg{hur.63}=2\incg{hur.64}
\end{equation}
with all possible colourings. There are three different ways to colour this
relation.

All three representations have the same dimension which we denote by $\delta$.
Then the challenge is to deduce directly from these relations that there is a finite set of possible values for $\delta$ which includes $\{0,1,2,4,8\}$.
The intention is not just to deduce this result but also to derive the additional relations
that hold in each case.

\section{Super Yang-Mills}
In this section we discuss an algebraic condition that only works in dimensions 3,4,6,10. A friendly discussion of this condition is given in \cite{YangMills}
(available as \href{http://arxiv.org/abs/0909.0551v2}{Division algebras and super symmetry}).

One way to distinguish these dimensions is that for any $V$ with a non-degenerate inner 
product we can construct a cubic Jordan algebra on $F\oplus V$. However the Jordan algebras
of 3,4,6,10 have the property that they are isomorphic to $H_2(\bK)$ for a composition algebra
$\bK$. This is discussed in \cite[\S 3.3]{MR1886087}. One would expect these cubic Jordan
algebras to satisfy an additional algebraic relation. There is a relation between these Jordan
algebras and the Clifford algebras which we discuss below. Namely the Jordan algebra is a
subalgebra of the Jordan algebra of the Clifford algebra.

In this section we stay with the same theme and move up another dimension.
Let $H$ be a two dimensional hyperbolic space. Then for any $V$ we have that
$\Cliff(V\oplus H)\cong M_2(\Cliff(V))$. Then this gives that $\Cliff(\bK\oplus H)$
acts on the direct sum of four copies of $\bK$ and $\Cliff^{\mathrm{even}}(\bK\oplus H)$
acts on $\bK\oplus \bK$. Equivalently $\bK\oplus \bK$
is a spin representation of $\Spin(\bK\oplus H)$. This is sometimes expressed as
\[ \Spin(\bK\oplus H)\cong \SL_2(\bK) \]
but this needs some interpretation for $\bK=\bH$ and even more interpretation for
$\bK=\bO$.

Fix a vector space $V$ with a non-degenerate inner product. Then a {\em system} consists
of a representation of $\Spin(V)$, $S$ and linear maps 
\[ \Gamma\colon S^*\otimes S^*\rightarrow V\qquad \Gamma\colon S\otimes S\rightarrow V \]
These maps are both required to be symmetric and to be maps of representations of $\Spin(V)$.
There is a further condition which in index notation is written
\[ \Gamma^\mu_{ab}\Gamma^{\nu bc}+\Gamma^\nu_{ab}\Gamma^{\mu bc} = 2g^{\mu\nu}\delta_a^c \]
Alternatively, there are linear maps
\[ \rho\colon V\otimes S^*\rightarrow S \qquad \rho\colon V\otimes S\rightarrow S^* \]
such that $\rho\oplus\rho\colon V\rightarrow \End(S^*\oplus S)$ satisfies the Clifford relations.
The equivalence between these definitions is given by
\begin{equation}\label{eq} \left\langle \rho(v\otimes s),t\right\rangle = (\Gamma(s\otimes t),v )
\end{equation}
for all $s,t\in S^*$ and all $v\in V$.

However there is an additional algebraic condition. This condition is needed to show that
the Lagrangian in supersymmetric Yang-Mills field theory is supersymmetric and is also 
needed in supersymmetric string theory.

This is the condition that the following quartic invariant vanishes identically
\[ s \mapsto |\Gamma(s,s)| \]
For further interpretations of this see \cite[(6.14),(6.69),(6.72)]{MR1701618}.
Polarising and taking into account the symmetries this is equivalent to
\begin{equation}\label{pol}
 (\Gamma(s,t),\Gamma(u,v))+ (\Gamma(s,u),\Gamma(t,v))+(\Gamma(s,v),\Gamma(t,u))=0
\end{equation} 
A system that satisfies this condition is {\em special}.

Next we construct an example from a composition algebra. Let $H$ be hyperbolic space
with basis $\{\alpha,\beta\}$ and quadratic form $\lambda\alpha+\mu\beta \mapsto \lambda\mu$.
Then the Clifford algebra is the algebra of $2\times 2$ matrices. An isomorphism is given by
\[ \alpha\mapsto \begin{pmatrix} 0&1\\ 0&0\end{pmatrix} \qquad
 \beta\mapsto \begin{pmatrix} 0&0\\ 1&0\end{pmatrix}
\]
Then by the theory of Clifford algebras there is a representation of
$\Cliff(\bK\oplus H)$ on $\bK\oplus\bK\oplus\bK\oplus\bK$. This is given explicitly
by
\begin{align*}
 \rho(a)\colon (u,v,x,y)&\mapsto (\overline{a}\,\overline{v},\overline{u}\,\overline{a},
-\overline{a}\,\overline{y},-\overline{x}\,\overline{a}) \\
\rho(\alpha)\colon (u,v,x,y)&\mapsto (x,y,0,0) \\
\rho(\beta)\colon (u,v,x,y)&\mapsto (0,0,u,v)
\end{align*}

The decomposition as a representation of $\Cliff^{\mathrm{even}}(V)$ is
\[ (u,v,x,y)=(u,y)+(v,x) \]
Then the definition of a system requires that these are dual representations.
A non-degenerate pairing is given by
\[ (u,y)\otimes (v,x) \mapsto -\left\langle u,\overline{x}\right\rangle +
\left\langle v,\overline{y}\right\rangle \]
 The proof that this pairing is a map of representations is left to the reader.

This shows that we have constructed a system. In the alternative formulation we have
a symmetric morphism of representations of $\Spin(\bK\oplus H)$,
\[ \Gamma\colon (\bK\oplus\bK)\otimes (\bK\oplus\bK)\rightarrow \bK\oplus H
\]

This is given by
\[ \Gamma\colon (u,y)\otimes (u^\prime,y^\prime)\mapsto (u^\prime y+y^\prime u) -
\left\langle u,u^\prime\right\rangle\alpha +\left\langle y,y^\prime\right\rangle\beta
\]

Then we check condition \eqref{eq} which determines $\Gamma$ uniquely.
This follows from the calculations
\begin{align*}
  \left\langle \rho(a\otimes(u,y)),(u^\prime,y^\prime) \right\rangle &=
\left\langle \overline{u}\,\overline{a},\overline{y}^\prime\right\rangle +
\left\langle u^\prime,ya\right\rangle \\
 \left\langle \rho(\alpha\otimes(u,y)),(u^\prime,y^\prime) \right\rangle 
&=\left\langle y,\overline{y}^\prime\right\rangle \\
 \left\langle \rho(\beta\otimes(u,y)),(u^\prime,y^\prime) \right\rangle 
&=-\left\langle u^\prime,\overline{u}\right\rangle
\end{align*}

Then we have $(u,y)\otimes (u,y) \mapsto uy+yu - |u|\alpha+|y|\beta$ and this has norm
$|uy|-|u|.|y|=0$. This shows that this system is special.

Next we show that there is a correspondence between trialities for $W$ and special systems
for $V=W\oplus H$. Assume we have a special system for $V=W\oplus H$.
Then we have a representation of $\Cliff(V)$.
There is a decomposition of $1$ into two orthogonal idempotents,
\[ 1 = \rho(\alpha)\rho(\beta)+\rho(\beta)\rho(\alpha) \]
These idempotents commute with $\Cliff(W)$ and so as a representation
of $\Cliff(W)$ this decomposes. Each of these can be decomposed
as a representation of $\Cliff^{\mathrm{even}}(W)$. This means the representation
can be written
\begin{align*}
 \rho(a)\colon (u,v,x,y)&\mapsto (\rho(a\otimes v),\overline{\rho}(a\otimes u),
-\rho(a\otimes y),-\overline{\rho}(a\otimes x) ) \\
\rho(\alpha)\colon (u,v,x,y)&\mapsto (x,y,0,0) \\
\rho(\beta)\colon (u,v,x,y)&\mapsto (0,0,u,v)
\end{align*}
The decomposition as a representation of $\Cliff^{\mathrm{even}}(V)$ is
\[ (u,v,x,y)=(u,y)+(v,x) \]
Then the definition of a system requires that these are dual representations.
A non-degenerate pairing is given by
\[ (u,y)\otimes (v,x) \mapsto -\left\langle u,x\right\rangle +
\left\langle v,y\right\rangle \]
This is required to be a morphism of representations. At this stage we have three vector
spaces $(W,S^+,S^-)$ each with a non-degenerate inner product and $S^+\oplus S^-$ is
a representation of $\Cliff(W)$.

Next we discuss $\Gamma$. This is given by
\[ \Gamma\colon (u,y)\otimes (u^\prime,y^\prime)\mapsto \gamma(u^\prime ,y)+ \gamma(y^\prime , u) -
\left\langle u,u^\prime\right\rangle\alpha +\left\langle y,y^\prime\right\rangle\beta
\]

Then the condition \eqref{eq} determines $\Gamma$. This is equivalent to the following
conditions which determine $\gamma$, $\overline{\gamma}$.
\begin{equation}
 \left\langle \rho(a\otimes u),v\right\rangle =
\left( a,\gamma(u\otimes v) \right) \qquad
\left\langle \overline{\rho}(a\otimes v),u\right\rangle =
\left( a,\overline{\gamma}(v\otimes u) \right) 
\end{equation}

This follows from the calculations
\begin{align*}
  \left\langle \rho(a\otimes(u,y)),(u^\prime,y^\prime) \right\rangle &=
\left\langle \overline{u}\,\overline{a},\overline{y}^\prime\right\rangle +
\left\langle u^\prime,ya\right\rangle \\
 \left\langle \rho(\alpha\otimes(u,y)),(u^\prime,y^\prime) \right\rangle 
&=\left\langle y,\overline{y}^\prime\right\rangle \\
 \left\langle \rho(\beta\otimes(u,y)),(u^\prime,y^\prime) \right\rangle 
&=-\left\langle u^\prime,\overline{u}\right\rangle
\end{align*}

Then it remains to show that the condition \eqref{pol} for $\Gamma$ is equivalent to
the Clifford relation for $\gamma\oplus\overline{\gamma}$. From the definition
\[ \Gamma( (u,y)\otimes (u,y) )= \gamma(u,y)+\overline{\gamma}(y,u)-|u|\alpha+|y|\beta \]
Then the condition that this has norm zero is equivalent to
\[ \left| \gamma(u,y)+\overline{\gamma}(y,u) \right| = |u|.|v| \]
Now consider $\gamma(u,y)\overline{\gamma}(y,u)+\overline{\gamma}(y,u)\gamma(u,y)$.
This is a morphism of representations $V\rightarrow V$ and $V$ is ireducible.
So by Schur's lemma this is a scalar matrix. Hence this is equivalent to 
condition \eqref{pol}.

Again we can express this in terms of diagrams. There are two types of edges.
An undirected edge for $V$ and a directed edge for $S$. Then there are two
trivalent vertices
\begin{equation}\label{gamma}
 \incg{hur.65}\qquad\incg{hur.66}
\end{equation}
Then the relations these satisfy are the Clifford relations and the relation \eqref{pol}.
The Clifford relations are
\begin{align*}
\incg{hur.67}+\incg{hur.68}&=2\incg{hur.69} \\
\incg{hur.70}+\incg{hur.71}&=2\incg{hur.72}
\end{align*}
The relation \eqref{pol} is
\begin{equation}\label{spc}
 \incg[0.5in]{hur.73}+\incg[0.5in]{hur.74}+\incg[0.5in]{hur.75}=0
\end{equation}

However these are not all the relations. There is also the requirement that $S$ is
a representation of $\Spin(V)$ and that these are morphisms of representations.
The way to understand these relations is to introduce a third type of edge for the
adjoint representation, $\fg$. Next we write the relations using this edge and then we
eliminate this edge using the fact that the adjoint representation is the exterior
square of $V$. Here we will not draw the diagrams that say that $\fg$ is a Lie algebra
and that $V$ is a representation as these relations are satisfied.

There are three trivalent vertices involving the new edge. These give the action of
$\fg$ on $V$, $S$, $S^*$, $\fg$. The four vertices are
\begin{equation}\label{ver}
 \incg{hur.76}\qquad\incg{hur.77}\qquad\incg{hur.78}\qquad\incg{hur.79}
\end{equation}
The first relations are the antisymmetry relations
\begin{equation}\label{as}
  \incg{hur.80}=-\incg{hur.81}\qquad \incg{hur.82}=-\incg{hur.83}
\end{equation}
This implies that $V$ is a self-dual representation and that
the pairing between $S$ and $S^*$ is a morphism of representations.

Next we have the relations that say that $V$, $S$, $S^*$, $\fg$ are representations.
This relation for $S$ is
\begin{equation}\label{ihx1}
 \incg{hur.84}=\incg{hur.85}+\incg{hur.86}
\end{equation}
The relation for $S^*$ is obtained by reversing the directions. This is a consequence
of \eqref{as} and \eqref{ihx1}.

Finally we have the relation that the maps \eqref{gamma} are morphisms of representations.
These are the relations
\begin{align*}\label{ihx1}
 \incg{hur.87}&=\incg{hur.88}+\incg{hur.89} \\
 \incg{hur.90}&=\incg{hur.91}+\incg{hur.92}
\end{align*}

Now we have obtained the relations the edge for the adjoint representation can be eliminated.
Three of the four vertices in \eqref{ver} are eliminated using
\begin{align*}
 \incg[0.5in]{hur.93}&\mapsto\frac12\incg[0.5in]{hur.94}-\frac12\incg[0.5in]{hur.95} \\
 \incg[0.5in]{hur.96}&\mapsto\frac12\incg[0.5in]{hur.97}-\frac12\incg[0.5in]{hur.98} \\
 \incg[0.5in]{hur.100}&\mapsto\frac12\incg[0.5in]{hur.101}-\frac12\incg[0.5in]{hur.102} 
\end{align*}

Again the challenge is to understand the consequences of these relations. It is reasonable
to expect that there is a finite set of possible values for $\delta$ which includes
$\{0,2,3,4,6,10\}$. As before the intention is not just to
prove this result but also to find additional relations which hold in each case.

If we omit the special relation \eqref{spc} then this can be realised for an infinite set
of values for $\delta$. Therefore this makes sense for $\delta$ an indeterminate.
Then from these relations we can formally construct a Lorentz supergroup for a spacetime
of dimension $\delta$. This is analogous to the way the Brauer category formally constructs
an orthogonal group for space of dimension $\delta$.

\newcommand{\etalchar}[1]{$^{#1}$}


\begin{thebibliography}{DEF{\etalchar{+}}99}

\bibitem[Bae02]{MR1886087}
John~C. Baez.
\newblock The octonions.
\newblock {\em Bull. Amer. Math. Soc. (N.S.)}, 39(2):145--205 (electronic),
  2002.

\bibitem[BH]{YangMills}
John~C. Baez and John Huerta.
\newblock Division algebras and supersymmetry i.
\newblock arXiv/0909.0551.

\bibitem[Boo98]{Diploma}
Dominik Boos.
\newblock {\em Ein tensorkategorieller Zugang zum Satz von Hurwitz}.
\newblock PhD thesis, ETH Zurich, 1998.

\bibitem[BS03]{MR2020553}
C.~H. Barton and A.~Sudbery.
\newblock Magic squares and matrix models of {L}ie algebras.
\newblock {\em Adv. Math.}, 180(2):596--647, 2003.

\bibitem[CS50]{MR0034378}
Claude Chevalley and R.~D. Schafer.
\newblock The exceptional simple {L}ie algebras {$F\sb 4$} and {$E\sb 6$}.
\newblock {\em Proc. Nat. Acad. Sci. U. S. A.}, 36:137--141, 1950.

\bibitem[Cvi08]{MR2418111}
Predrag Cvitanovi{\'c}.
\newblock {\em Group theory}.
\newblock Princeton University Press, Princeton, NJ, 2008.
\newblock Birdtracks, Lie's, and exceptional groups.

\bibitem[DEF{\etalchar{+}}99]{MR1701618}
Pierre Deligne, Pavel Etingof, Daniel~S. Freed, Lisa~C. Jeffrey, David Kazhdan,
  John~W. Morgan, David~R. Morrison, and Edward Witten, editors.
\newblock {\em Quantum fields and strings: a course for mathematicians. {V}ol.
  1, 2}.
\newblock American Mathematical Society, Providence, RI, 1999.
\newblock Material from the Special Year on Quantum Field Theory held at the
  Institute for Advanced Study, Princeton, NJ, 1996--1997.

\bibitem[DP08]{MR2492402}
Carlos Dur{\'a}n and Thomas P{\"u}ttmann.
\newblock A minimal {B}rieskorn 5-sphere in the {G}romoll-{M}eyer sphere and
  its applications.
\newblock {\em Michigan Math. J.}, 56(2):419--451, 2008.

\bibitem[Hur98]{hurwitz}
Adolf Hurwitz.
\newblock Über die composition der quadratischen formen von beliebig vielen
  variabeln.
\newblock {\em Nachr. Ges. Wiss. Göttingen}, 1898.

\bibitem[Kup94]{MR1265145}
Greg Kuperberg.
\newblock The quantum {$G\sb 2$} link invariant.
\newblock {\em Internat. J. Math.}, 5(1):61--85, 1994.

\bibitem[Kup96]{MR1403861}
Greg Kuperberg.
\newblock Spiders for rank {$2$} {L}ie algebras.
\newblock {\em Comm. Math. Phys.}, 180(1):109--151, 1996.

\bibitem[Mas83]{MR723943}
W.~S. Massey.
\newblock Cross products of vectors in higher-dimensional {E}uclidean spaces.
\newblock {\em Amer. Math. Monthly}, 90(10):697--701, 1983.

\bibitem[Ros96]{MR1397790}
Markus Rost.
\newblock On the dimension of a composition algebra.
\newblock {\em Doc. Math.}, 1:No. 10, 209--214 (electronic), 1996.

\bibitem[Tit66]{MR0214638}
J.~Tits.
\newblock Sur les constantes de structure et le th\'eor\`eme d'existence des
  alg\`ebres de {L}ie semi-simples.
\newblock {\em Inst. Hautes \'Etudes Sci. Publ. Math.}, (31):21--58, 1966.

\bibitem[Wes07]{MR2320368}
Bruce~W. Westbury.
\newblock Enumeration of non-positive planar trivalent graphs.
\newblock {\em J. Algebraic Combin.}, 25(4):357--373, 2007.

\end{thebibliography}
\end{document}